\pdfoutput=1
\documentclass{article}
\usepackage[utf8]{inputenc}
\usepackage{amsmath}
\usepackage{amssymb}
\usepackage{amsthm}
\usepackage{enumerate}
\usepackage{mathrsfs}
\usepackage{scrextend}
\usepackage{mathtools}
\usepackage[hyperfootnotes=false]{hyperref}
\usepackage{cleveref}
\usepackage[multiple]{footmisc}
\usepackage[a4paper, total={6in, 9in}]{geometry}
\usepackage{scalerel,stackengine}
\stackMath
\newcommand\reallywidehat[1]{%
\savestack{\tmpbox}{\stretchto{%
  \scaleto{%
    \scalerel*[\widthof{\ensuremath{#1}}]{\kern-.6pt\bigwedge\kern-.6pt}%
    {\rule[-\textheight/2]{1ex}{\textheight}}
  }{\textheight}%
}{0.5ex}}%
\stackon[1pt]{#1}{\tmpbox}%
}
\theoremstyle{plain}
\newtheorem{theorem}{Theorem}
\newtheorem{lemma}[theorem]{Lemma}
\newtheorem{proposition}[theorem]{Proposition}

\newtheorem{problem}[theorem]{Problem}

\theoremstyle{definition}

\newtheorem{question}[theorem]{Question}

\numberwithin{theorem}{section}

\title{On a problem of Erdős and Graham about consecutive sums in strictly increasing sequences}
\author{Adrian Beker\footnote{University of Zagreb, Faculty of Science, Department of Mathematics, Zagreb,
Croatia.\\ Email: \nolinkurl{adrian.beker@math.hr}}}
\date{\today}

\begin{document}

\maketitle

\begin{abstract}
    We show the existence of a constant $c > 0$ such that, for all positive integers $n$, there exist integers $1 \leq a_1 < \ldots < a_k \leq n$ such that there are at least $cn^2$ distinct integers of the form $\sum_{i=u}^{v}a_i$ with $1 \leq u \leq v \leq k$. This answers a question of Erdős and Graham. We also prove a non-trivial upper bound on the maximum number of distinct integers of this form and address several open problems.
\end{abstract}

\section{Introduction}

Given a finite sequence of integers $a = (a_i)_{1\leq i \leq k}$, we denote by $S(a)$ the set of its consecutive sums, that is,
$$S(a) \vcentcolon= \Biggl\{\sum_{i=u}^{v}a_i\ \Bigg|\ 1 \leq u \leq v \leq k\Biggr\}.$$
By relating $S(a)$ to the difference set of the set of partial sums of $a$, one immediately notices that it is a natural object to study from an additive-combinatorial perspective. The study of $S(a)$ in various settings was initiated by Erdős, who together with Graham posed the following problem in \cite{erdos-graham}. This problem also appears as problem \#356 on Thomas Bloom's website \emph{Erdős problems} \cite{bloom}.

\begin{problem}
\label{main problem}
Is there some $c > 0$ such that, for all sufficiently large $n$, there exist integers $1 \leq a_1 < \ldots < a_k \leq n$ such that there are at least $cn^2$ distinct integers of the form $\sum_{i=u}^{v}a_i$ with $1 \leq u \leq v \leq k$?
\end{problem}

The obvious example in which $k = n$ and $a_i = i$ for all $1 \leq i \leq k$ only just fails, attaining $\Theta(n^2(\log n)^{-\delta+o(1)})$ distinct consecutive sums, where $\delta = 1-\frac{1+\log\log2}{\log2} \approx 0.086$ is the Erdős-Ford-Tenenbaum constant. This is a straightforward consequence of Ford's bounds on the multiplication table problem (see \cite{ford1} and \cite{ford2}). Erdős and Graham also asked what happens if we drop the monotonicity assumption and whether the same is true for permutations of $[n]$ instead of increasing sequences. In the former case, an affirmative answer was given by Hegyvári, who showed more strongly in \cite{hegyvari} that one can find a sequence of length $\left(\frac{1}{3}+o(1)\right)n$ in $[n]$ with all consecutive sums distinct. In the latter case, Konieczny \cite{konieczny} estalished the rather striking fact that, if $a$ is a uniform random permutation of $[n]$, then $|S(a)| \sim \Bigl(\frac{1+e^{-2}}{4}\Bigr)n^2$ with high probability.

However, despite previous related work, Problem \ref{main problem} remained open. In this paper, we solve this problem by proving the following result:

\begin{theorem}
\label{main result}
There exists a constant $c_1 > 0$ such that, for all positive integers $n$, there exist integers $1 \leq a_1 < \ldots < a_k \leq n$ such that there are at least $c_1n^2$ distinct integers of the form $\sum_{i=u}^{v}a_i$ with $1 \leq u \leq v \leq k$.
\end{theorem}

We establish Theorem \ref{main result} via the following probabilistic result, thereby showing that, in a sense, there are many sequences $a$ exhibiting the desired behaviour.

\begin{theorem}
\label{probabilistic version of main result}
There exists a constant $c_2 > 0$ such that the following holds for all positive integers $n$. Let $\varepsilon_1, \ldots, \varepsilon_n$ be i.i.d.\ Rademacher random variables and set $a_i = 3i + \varepsilon_i$ for $1 \leq i \leq n$. Then with positive probability, we have $|S(a)| \geq c_2n^2$.
\end{theorem}

We also provide a class of explicit examples of strictly increasing sequences with many consecutive sums:

\begin{theorem}
\label{deterministic version of main result}
There exists a constant $c_3 > 0$ such that the following holds for all positive integers $n$. Let $b$ be a positive integer such that $\log n \leq b \leq \frac{n}{(\log n)^2}$ and define
$$a_i = \begin{cases}2i & \text{if } b \mid i\\2i-1 & \text{otherwise}\end{cases}$$
for $1 \leq i \leq n$. Then $|S(a)| \geq c_3n^2$.
\end{theorem}

We prove Theorems \ref{probabilistic version of main result} and \ref{deterministic version of main result} by showing that the additive energy of the set of partial sums of $a$ is essentially as small as possible (up to a multiplicative constant). Roughly speaking, the idea is that if there aren't many pairs of consecutive sums that are equal, then there must be many distinct consecutive sums. We give a full proof of Theorem \ref{probabilistic version of main result} and then sketch the modifications that are needed to obtain Theorem \ref{deterministic version of main result}. These arguments are carried out in Section 2. 

Finally, in Section 3, we make some concluding remarks about our work and discuss several questions that remain open. In particular, we show that the trivial upper-bound $\frac{n(n+1)}{2}$ on $|S(a)|$ is not sharp:

\begin{proposition}
\label{upper bound on number of consecutive sums}
Let $n$ be a positive integer and let $1 \leq a_1 < \ldots < a_k \leq n$ be integers. Then 
$$|S(a)| \leq (c_4+o(1))n^2,$$
where $c_4 = \frac{e^2-1}{2(e^2+1)} \approx 0.381$.
\end{proposition}

\noindent\textbf{Notation.} We use standard asymptotic notation. Given functions $f, g \colon \mathbb{N} \to \mathbb{C}$, we write $f = O(g)$ or $f \ll g$ if there is a constant $C > 0$ such that $|f(n)| \leq C|g(n)|$ for all sufficiently large $n$. We also write $f = o(g)$ to mean that $\frac{f(n)}{g(n)} \to 0$ as $n \to \infty$. If $f = O(g)$ and $g = O(f)$, then we write $f = \Theta(g)$.

For a positive integer $n$, we abbreviate the set $\{1,\ldots,n\}$ to $[n]$. For sets $A, B \subseteq \mathbb{Z}$, we define their difference set to be $A - B \vcentcolon= \{a - b \mid a \in A,\ b \in B\}$. Given $a \in \mathbb{Z}$ and $m \in \mathbb{N}$, we write $[a \text{ mod } m]$ for the unique integer $r$ such that $0 \leq r \leq m-1$ and $a \equiv r \pmod m$. For the sake of simplicity, we abuse the notation $[a, b]$ to mean $[a,b] \cap \mathbb{Z}$ and abuse the word \emph{interval} by calling such sets intervals.

\section{Proof of Theorems \ref{probabilistic version of main result} and \ref{deterministic version of main result}}

The key to the proof of Theorems \ref{probabilistic version of main result} and \ref{deterministic version of main result} is the notion of \emph{additive energy} and its relation to the size of sumsets/difference sets. Following Tao and Vu \cite{tao-vu}, we define the additive energy of a finite non-empty set $P \subseteq \mathbb{Z}$ to be
$$E(P) \vcentcolon= |\{(x,y,z,w) \in P^4 \mid x-y = z-w\}|.$$
Writing $r_P(t)$ for the number of representations of $t \in \mathbb{Z}$ as a difference of two elements of $P$, one obtains the following expression for the additive energy:
$$E(P) = \sum_{t\in P-P}r_P(t)^2.$$
Since $\sum_{t\in P-P}r_P(t) = |P|^2$, the Cauchy-Schwarz inequality implies that 
$$E(P) \geq \frac{|P|^4}{|P-P|}.$$
Thus, if the difference set is small, the additive energy must be large. Contrapositively, $P-P$ must be large in order for $E(P)$ to be small. In our setting, given a finite sequence of positive integers $a = (a_1, \ldots, a_n)$, we define the corresponding sequence of partial sums $(p_0, p_1, \ldots, p_n)$ by $p_i = \sum_{j=1}^{i}a_j$ for $0 \leq i \leq n$. In particular, $0 = p_0 < p_1 < \ldots < p_n$. Thus, on writing $P(a) = \{p_0, p_1, \ldots, p_n\}$, we have
$$S(a) = \{p_j - p_i \mid 0 \leq i < j \leq n\} = (P(a)-P(a)) \cap \mathbb{N},$$ 
so in particular $|P(a)| = n+1$ and $|S(a)| = \frac{|P(a)-P(a)|-1}{2}$. Therefore, in order to establish Theorem \ref{probabilistic version of main result}, it suffices to prove the following:

\begin{theorem}
\label{expected additive energy}
Let $n$ be a positive integer and let $\varepsilon_1, \ldots, \varepsilon_n$ be i.i.d.\ Rademacher random variables. Define $a_i = 3i + \varepsilon_i$ for $1 \leq i \leq n$. Then the expected value of $E(P(a))$ is $O(n^2)$.
\end{theorem}

Before going on to prove Theorem \ref{expected additive energy}, we establish a simple bound on the probability of a symmetric binomial random variable being divisible by a given integer.

\begin{lemma}
\label{facts about the binomial distribution}
Let $m, n$ be positive integers and let $X$ be a binomial random variable with parameters $n$ and $\frac{1}{2}$. Then $\mathbb{P}(X \equiv 0\ (\mathrm{mod}\ m)) \leq \frac{1}{m} + \frac{2}{\sqrt{n}}$.
\end{lemma}
\noindent\textit{Proof.} We have
$$\mathbb{P}(X \equiv 0\ (\mathrm{mod}\ m)) = \sum_{\substack{0\leq k\leq n\\ k \equiv 0\ (\mathrm{mod}\ m)}}f(k),$$
where $f$ is the probability mass function of $X$. We split this sum as follows:
$$\mathbb{P}(X \equiv 0\ (\mathrm{mod}\ m)) = \sum_{\substack{0\leq k\leq \frac{n}{2}-m\\k \equiv 0\ (\mathrm{mod}\ m)}}f(k) + \sum_{\substack{\frac{n}{2}+m\leq k \leq n\\k \equiv 0\ (\mathrm{mod}\ m)}}f(k) + \sum_{\substack{\left|k-\frac{n}{2}\right|<m\\k \equiv 0\ (\mathrm{mod}\ m)}}f(k).$$
Since $f$ is increasing on $\left[0,\frac{n}{2}\right]$ and decreasing on $\left[\frac{n}{2},n\right]$, the first two sums can be upper bounded as follows:
$$\sum_{\substack{0\leq k\leq \frac{n}{2}-m\\k \equiv 0\ (\mathrm{mod}\ m)}}f(k) \leq \frac{1}{m}\sum_{\substack{0\leq k\leq \frac{n}{2}-m\\k \equiv 0\ (\mathrm{mod}\ m)}}\sum_{l=k}^{k+m-1}f(l) \leq \frac{1}{m}\sum_{0\leq l \leq \frac{n}{2}-1}f(l),$$
$$\sum_{\substack{\frac{n}{2}+m\leq k \leq n\\k \equiv 0\ (\mathrm{mod}\ m)}}f(k) \leq \frac{1}{m}\sum_{\substack{\frac{n}{2}+m\leq k \leq n\\k \equiv 0\ (\mathrm{mod}\ m)}}\sum_{l=k-m+1}^{k}f(l) \leq \frac{1}{m}\sum_{\frac{n}{2}+1 \leq 
l \leq n}f(l).$$
On the other hand, by standard estimates on central binomial coefficients, the maximum of $f$ is at most $\frac{1}{\sqrt{n}}$. Since there are at most two multiples of $m$ strictly between $\frac{n}{2}-m$ and $\frac{n}{2}+m$, we have
$$\sum_{\substack{\left|k-\frac{n}{2}\right|<m\\k \equiv 0\ (\mathrm{mod}\ m)}}f(k) \leq \frac{2}{\sqrt{n}}.$$
The conclusion follows by putting together the obtained estimates. $\qed$
\\\\
\noindent\textbf{Remark.} The quantity of interest in Lemma \ref{facts about the binomial distribution} is closely related to the $n$-step transition probabilities of a simple symmetric random walk on $\mathbb{Z}/m\mathbb{Z}$. By using Fourier analysis instead of ad-hoc arguments, one can obtain similar bounds for a general symmetric step distribution with bounded support. Since the rest of the proof of Theorem \ref{expected additive energy} carries over verbatim to this setting, one obtains a more general version of Theorem \ref{probabilistic version of main result}.
\\\\
\noindent\textit{Proof of Theorem \ref{expected additive energy}.} We begin by deriving a more convenient expression for the additive energy of $P(a)$:
\begin{align*}
    E(P(a)) &= |\{(i,j,k,l) \in [0,n]^4 \mid p_j-p_i = p_l-p_k\}|\\
    &= 2|\{(i,j,k,l) \in [0,n]^4 \mid i < j,\ k < l,\ p_j-p_i = p_l-p_k\}| + (n+1)^2\\
    &= 2\Biggl|\Biggl\{(i,j,k,l) \in [0,n]^4\ \Bigg|\ i < j,\ k < l,\ \sum_{u=i+1}^{j}a_u = \sum_{v=k+1}^{l}a_v\Biggr\}\Biggr| + (n+1)^2.
\end{align*}
By linearity of expectation, our task therefore reduces to showing that
$$\sum_{\substack{i,j,k,l \in [0,n]\\ i < j,\ k < l}}\mathbb{P}\Biggl(\sum_{u=i+1}^{j}a_u = \sum_{v=k+1}^{l}a_v\Biggr) = O(n^2).$$
We now make a further simplification by observing that, in the above sum, we may discard the pairs of intervals $[i+1, j]$, $[k+1, l]$ that intersect. Indeed, by symmetry, we may assume that $i \leq k$ (the $O$-notation takes care of the fact that the sum doubles if we include the pairs with $i > k$). Furthermore, we may ignore the terms in which either the intervals $[i+1, j]$, $[k+1, l]$ are equal or one of them is strictly contained in the other. Indeed, in the former case, there are $\frac{n(n+1)}{2}$ such terms and each of them contributes $1$ to the sum, whereas in the latter case, all terms are zero. This leaves us with a sum over all $i,j,k,l \in [0,n]$ such that either $i < k < j < l$ or $i < j \leq k < l$. But note that if $i < k < j < l$, then the equality 
$$\sum_{u=i+1}^{j}a_u = \sum_{v=k+1}^{l}a_v$$ can be rewritten as
$$\sum_{u=i+1}^{k}a_u = \sum_{v=j+1}^{l}a_v.$$ Hence, if we keep only the terms with $i < j \leq k < l$, our sum decreases by a factor of at most $2$. It follows that it is enough to show that
$$\sum_{\substack{i,j,k,l \in [0,n]\\ i < j \leq k < l}}\mathbb{P}\Biggl(\sum_{u=i+1}^{j}a_u = \sum_{v=k+1}^{l}a_v\Biggr) = O(n^2).$$
To this end, observe that
$$\mathbb{P}\Biggl(\sum_{u=i+1}^{j}a_u = \sum_{v=k+1}^{l}a_v\Biggr) = \mathbb{P}\Biggl(\sum_{u=i+1}^{j}\varepsilon_u-\sum_{v=k+1}^{l}\varepsilon_v = 3\Biggl(\sum_{v=k+1}^{l}v-\sum_{u=i+1}^{j}u\Biggr)\Biggr),$$
which, since the intervals $[i+1, j]$, $[k+1, l]$ are crucially disjoint, equals 
$$g\Biggl(|[i+1, j]| + |[k+1, l]|, 3\Biggl(\sum_{v\in[k+1, l]}v-\sum_{u\in[i+1, j]}u\Biggr)\Biggl).$$
Here, we let $g(m,\cdot)$ denote the probability mass function of the sum of $m$ i.i.d.\ Rademacher random variables. Hence, it is certainly enough to show that the sum of $g\Bigl(|I|+|J|, 3\big(\sum_{v\in J}v-\sum_{u\in I}u\big)\Bigr)$ over all pairs of intervals $I,J$ in $[1,n]$ is $O(n^2)$. By restricting to intervals of fixed lengths $k,l \in [1,n]$, the sum in question becomes
$$\sum_{i=1}^{n-k+1}\sum_{j=1}^{n-l+1}g\Biggl(k+l,3\Biggl(ki+\frac{k(k-1)}{2}-lj-\frac{l(l-1)}{2}\Biggr)\Biggr),$$
or in other words
$$\sum_{i=1}^{n-k+1}\sum_{j=1}^{n-l+1}f\Bigl(k+l,\frac{3}{2}(ik-jl)+\frac{3}{4}k(k-1)-\frac{3}{4}l(l-1)+\frac{1}{2}(k+l)\Bigr),$$
where $f(m,\cdot)$ denotes the probability mass function of the binomial distribution with parameters $m$ and $\frac{1}{2}$. By grouping equal terms together, this can be further rewritten as
\begin{equation}\label{inner product of number of representations with mass function}
    \sum_{x=0}^{k+l}p(x)f(k+l,x),
\end{equation}
where $p(x)$ denotes the number of pairs $(i,j) \in [1,n-k+1]\times[1,n-l+1]$ such that 
\begin{equation}\label{linear diophantine}
    ik-jl = \frac{1}{3}(2x-k-l)-\frac{1}{2}k(k-1)+\frac{1}{2}l(l-1).
\end{equation}
We now proceed to bound the sum (\ref{inner product of number of representations with mass function}) by analysing the support and the maximum of $p$. Writing $q \vcentcolon= \gcd(k,l)$, we see that for the linear Diophantine equation (\ref{linear diophantine}) to have at least one solution, it is necessary that $q$ divides $4x$, that is, $x$ is a multiple of $q' \vcentcolon = \frac{q}{\gcd(q,4)}$. On the other hand, if $p(x) > 0$, then the equation (\ref{linear diophantine}) can be rewritten as $ik'-jl' = t$, where $k' \vcentcolon= \frac{k}{q}$, $l' \vcentcolon= \frac{l}{q}$ and $t$ is some integer that depends on $x$. Since $k'$ is invertible modulo $l'$, all admissible values of $i$ come from a single congruence class modulo $l'$, so it follows that
$$p(x) \leq \left\lceil\frac{n-k+1}{l'}\right\rceil \leq \frac{n}{l'}+1 = \frac{qn}{l}+1.$$
Therefore, the sum (\ref{inner product of number of representations with mass function}) can be upper bounded as follows:
\begin{align*}
    \sum_{x=0}^{k+l}p(x)f(k+l,x) \leq \sum_{\substack{x=0\\x \equiv 0\ (\mathrm{mod}\ q')}}^{k+l}\left(\frac{qn}{l}+1\right)f(k+l,x)\leq 1 + \frac{n}{l}\left(4+\frac{2q}{\sqrt{k+l}}\right),
\end{align*}
where we used Lemma \ref{facts about the binomial distribution} in the second inequality. Finally, by interchanging the roles of $k$ and $l$ if necessary, we see that it suffices to show that
$$\sum_{1\leq k \leq l \leq n}\left[1+\frac{n}{l}\left(1+\frac{\gcd(k,l)}{\sqrt{k+l}}\right)\right] = O(n^2).$$
Since $\sum_{1\leq k \leq l \leq n}1 = \frac{n(n+1)}{2}$ and 
$$\sum_{1\leq k\leq l \leq n}\frac{n}{l} = \sum_{l=1}^{n}\sum_{k=1}^{l}\frac{n}{l} = \sum_{l=1}^{n}l\cdot \frac{n}{l} = n^2,$$
this reduces to showing that
$$\sum_{1\leq k \leq l\leq n}\frac{\gcd(k,l)}{l^{3/2}} = O(n),$$
which turns out to be an elementary manipulation involving arithmetic functions. To begin, observe that
\begin{align*}
    \sum_{l=1}^{n}\sum_{k=1}^{l}\frac{\gcd(k,l)}{l^{3/2}} = \sum_{l=1}^{n}\frac{1}{l^{3/2}}\sum_{d\mid l}d\varphi\left(\frac{l}{d}\right) = \sum_{l=1}^{n}\frac{1}{\sqrt{l}}\sum_{d\mid l}\frac{\varphi(l/d)}{l/d} = \sum_{l=1}^{n}\frac{1}{\sqrt{l}}\sum_{d'\mid l}\frac{\varphi(d')}{d'},
\end{align*}
where $\varphi$ is Euler's totient function. Hence, by interchanging the order of summation, we obtain that this equals
\begin{align*}
    \sum_{d=1}^{n}\frac{\varphi(d)}{d}\sum_{\substack{l=1\\ d \mid l}}^{n}\frac{1}{\sqrt{l}} = \sum_{d=1}^{n}\frac{\varphi(d)}{d^{3/2}}\sum_{l'=1}^{\lfloor n/d\rfloor}\frac{1}{\sqrt{l'}}.
\end{align*}
By employing the standard estimate 
$$\sum_{l'=1}^{L}\frac{1}{\sqrt{l'}} \leq \int_{0}^{L}\frac{1}{\sqrt{x}}\,dx = 2\sqrt{L}$$ 
and the trivial bound $\varphi(d) \leq d$, this can be upper bounded by
\begin{align*}
\sum_{d=1}^{n}\frac{\varphi(d)}{d^{3/2}}\cdot 2\sqrt{\frac{n}{d}} = 2\sqrt{n}\sum_{d=1}^{n}\frac{\varphi(d)}{d^2} \leq 2\sqrt{n}\sum_{d=1}^{n}\frac{1}{d} = O(\sqrt{n}\log n),
\end{align*}
so we are done. $\qed$

We now turn our attention to Theorem \ref{deterministic version of main result}. The main idea of the proof is similar as in the case of Theorem \ref{probabilistic version of main result}, so we will be fairly brief on the details. As before, the idea is to show that the additive energy of $P(a)$ is $O(n^2)$. In this case, we have the following general formula for consecutive sums:
$$\sum_{i=u+1}^{v}a_i = v^2 - u^2 + \left\lfloor\frac{v-u}{b}\right\rfloor + \varepsilon,$$ 
for some $\varepsilon \in \{0,1\}$ (depending on $u,v$). Thus, again fixing the lengths of the two intervals in consideration to be $k, l \in [n]$, our task amounts to bounding the number of pairs $(i,j) \in [n-k+1] \times [n-l+1]$ for which
\begin{equation}\label{another diophantine}
    2ki - 2lj = l^2 - k^2 + \left\lfloor\frac{l}{b}\right\rfloor - \left\lfloor\frac{k}{b}\right\rfloor + \delta
\end{equation}
for some $\delta \in \{-1,0,1\}$. Again letting $q \vcentcolon= \gcd(k, l)$ and writing $k' \vcentcolon= \frac{k}{q}$, $l' \vcentcolon= \frac{l}{q}$, we see that for (\ref{another diophantine}) to have at least one integer solution $(i,j)$, it is necessary that $\left\lfloor\frac{k}{b}\right\rfloor - \left\lfloor\frac{l}{b}\right\rfloor \equiv \delta \pmod q$ for some $\delta \in \{-1,0,1\}$. After multiplying through by $b$ and using that $b\left\lfloor\frac{m}{b}\right\rfloor = m - [m \text{ mod } b]$ for $m \in \mathbb{Z}$, this condition becomes equivalent to $[k \text{ mod } b] - [l \text{ mod } b] + \delta b$ being a multiple of $q$ and $k'-l'$ at the same time being congruent to $\frac{[k \text{ mod } b] - [l \text{ mod } b] + \delta b}{q}$ modulo $b$. In particular, fixing $q$, only the pairs $(k,l)$ with $k'-l' \in \Bigl(-\frac{2b}{q},\frac{2b}{q}\Bigr)+b\mathbb{Z}$ contribute to the additive energy. Therefore, it suffices to show that
$$\sum_{q=1}^{n}\sum_{\substack{1 \leq k' \leq l' \leq \frac{n}{q}\\k'-l' \in (-2b/q,2b/q)+b\mathbb{Z}}}\left(\frac{n}{l'}+1\right) = O(n^2).$$
Since for each $l'$ there are at most $\left(\frac{l'}{b}+1\right)\left(\frac{4b}{q}+1\right)$ values $k' \leq l'$ satisfying the congruence condition, this sum can be upper bounded by
\begin{align*}
    \sum_{q=1}^{n}\sum_{1\leq l'\leq \frac{n}{q}}\left(\frac{l'}{b}+1\right)\left(\frac{4b}{q}+1\right)\left(\frac{n}{l'}+1\right) &\ll n\sum_{q=1}^{n}\sum_{1\leq l'\leq \frac{n}{q}}\left(\frac{1}{b}+\frac{1}{l'}\right)\left(\frac{4b}{q}+1\right)\\
    &\ll n\sum_{\substack{q,l'\in[n]\\ql' \leq n}}\left(\frac{1}{q}+\frac{1}{b}+\frac{b}{ql'}+\frac{1}{l'}\right).
\end{align*}
But note that
$$\sum_{\substack{q,l'\in[n]\\ql' \leq n}}\frac{1}{q} = \sum_{\substack{q,l'\in[n]\\ql' \leq n}}\frac{1}{l'} = \sum_{a=1}^{n}\left\lfloor\frac{n}{a}\right\rfloor\cdot\frac{1}{a} \leq n\sum_{a=1}^{n}\frac{1}{a^2} = O(n),$$
whereas
$$\sum_{\substack{q,l'\in[n]\\ql' \leq n}}\frac{1}{b} = \frac{1}{b}\sum_{a=1}^{n}\left\lfloor\frac{n}{a}\right\rfloor \ll \frac{n\log n}{b}.$$
Finally, we have
$$\sum_{\substack{q,l'\in[n]\\ql' \leq n}}\frac{b}{ql'} \leq b\Biggl(\sum_{q=1}^{n}\frac{1}{q}\Biggr)\Biggl(\sum_{l'=1}^{n}\frac{1}{l'}\Biggr) \ll b(\log n)^2,$$
so Theorem \ref{deterministic version of main result} follows on combining the obtained estimates.

\section{Concluding remarks and open problems}

By establishing Theorem \ref{main result}, we have made progress on the problem of estimating the maximum of $|S(a)|$ over all strictly increasing sequences $(a_i)_{1\leq i \leq k}$ in $[n]$. A rough calculation shows that one may take $c_2 = 2\cdot10^{-2}$ in Theorem \ref{probabilistic version of main result} and hence $c_3 = 2\cdot 10^{-3}$ in Theorem \ref{main result}. This can certainly be improved upon by performing more careful calculations, but we spend no effort in doing so. On the other hand, we can complement this lower bound by proving a non-trivial upper bound in the form of Proposition \ref{upper bound on number of consecutive sums}. The argument presented here is reminiscent of the proofs of the upper bound in Theorem 1.2 and of Proposition 5.1 in \cite{konieczny}, albeit significantly simpler since it does not require any prior preparation.

\noindent\textit{Proof of Proposition \ref{upper bound on number of consecutive sums}.} Fix a parameter $\alpha \in (0,1)$, the exact value of which will be determined later. We split $S(a)$ into two parts: the elements that are less than $\frac{1}{2}\alpha(n+1)^2$ and those greater than or equal to $\frac{1}{2}\alpha(n+1)^2$. The former part has cardinality at most $\frac{1}{2}\alpha(n+1)^2$. The cardinality of the latter part does not exceed the number of pairs $(i,j) \in [k]^2$ with $\sum_{u=i}^{j}a_u \geq \frac{1}{2}\alpha(n+1)^2$. This number, in turn, is at most $|L_n|$, where we define
$$L_n \vcentcolon= \Biggl\{(i,j)\ \Bigg|\ 0 \leq i < j \leq n,\ \sum_{u=i+1}^{j}u \geq \frac{1}{2}\alpha(n+1)^2 \Biggr\}.$$
But we can rewrite $L_n$ as
$$L_n = \Bigg\{(i,j)\ \Bigg|\ 0 \leq i < j \leq n,\ \Bigl(\frac{j+1/2}{n+1}\Bigr)^2 - \Bigl(\frac{i+1/2}{n+1}\Bigr)^2 \geq \alpha\Bigg\}.$$
Thus, $L_n$ naturally corresponds to the set of points of the lattice $\Bigl(\frac{1}{n+1}(\mathbb{Z}+\frac{1}{2})\Bigr)^2$ inside the set $\Lambda_{\alpha}$, where we define
$$\Lambda_{\alpha} \vcentcolon= \Big\{(x,y) \in [0,1]^2\ \Big|\ y^2 - x^2 \geq \alpha\Big\}.$$
On associating to each such point $(x,y)$ the square $(x-\frac{1}{n+1},x]\times(y,y+\frac{1}{n+1}]$, a standard volume packing argument shows that
$$\frac{|L_n|}{(n+1)^2} \leq |\Lambda_\alpha| + O\Bigl(\frac{1}{n+1}\Bigr).$$
Here, if $E \subseteq \mathbb{R}^2$ is a measurable set, we denote by $|E|$ its Lebesgue measure. Hence, we have $|L_n| \leq (|\Lambda_{\alpha}|+o(1))(n+1)^2$, and a simple calculation shows that
$$|\Lambda_{\alpha}| = \int_{0}^{\sqrt{1-\alpha}}\Bigl(1-\sqrt{x^2+\alpha}\Bigr)\,dx = \frac{1}{2}\Biggl(\sqrt{1-\alpha}-\alpha\log\Biggl(\frac{1+\sqrt{1-\alpha}}{\sqrt{\alpha}}\Biggr)\Biggr).$$
Therefore, we may take $c_4 = \frac{1}{2}h(\alpha)$, where
$$h \colon (0,1) \to \mathbb{R}, \quad \alpha \mapsto \alpha+\sqrt{1-\alpha}-\alpha\log\Biggl(\frac{1+\sqrt{1-\alpha}}{\sqrt{\alpha}}\Biggr).$$
All that remains is to optimise the function $h$, which can be done by routine calculus. One finds that the derivative of $h$ is $h'(\alpha) = 1 - \log\Bigl(\frac{1+\sqrt{1-\alpha}}{\sqrt{\alpha}}\Bigr)$, so $h$ has a minimum at $\alpha = \left(\frac{2e}{e^2+1}\right)^2$. Therefore, the minimum value of $h$ is $\frac{e^2-1}{e^2+1}$, and the conclusion follows. $\qed$

In spite of Theorem \ref{main result} and Proposition \ref{upper bound on number of consecutive sums}, the upper and lower bounds on the maximum of $|S(a)|$ remain quite far apart. Hence, we ask the following question.

\begin{question}
\label{asymptotics for maximum number of consecutive sums}
Let $\mathcal{A}_n$ denote the family of strictly increasing sequences in $[n]$. What is $\max_{a \in \mathcal{A}_n}|S(a)|$? In particular, does there exist a constant $c > 0$ such that 
$$\max_{a \in \mathcal{A}_n}|S(a)| = (c+o(1))n^2,$$ 
and if so, what is the value of $c$?
\end{question}

Perhaps even more interesting is the question of determining the extent to which the property of having many consecutive sums is typical of strictly increasing sequences in $[n]$ of linear length. Theorem \ref{probabilistic version of main result} suggests that this property might be typical. It would be interesting to investigate whether this is true for sequences arising from the binomial model for random dense subsets of $[n]$. Specifically, given a parameter $p \in [0,1]$, a \emph{$p$-random subset} of $[n]$ is a random set obtained by taking each element of $[n]$ independently at random with probability $p$.

\begin{question}
\label{consecutive sums of random sequences}
Is it true that, for any $p \in (0, 1)$, there exists a constant $c > 0$ such that, if $A = \{a_1, \ldots, a_k\}$ is a $p$-random subset of $[n]$ with $a_1 < \ldots < a_k$, then $|S(a)| \geq cn^2$ with high probability as $n \to \infty$?
\end{question}

Before making an attempt to answer Question \ref{consecutive sums of random sequences}, it might be instructive to consider the corresponding question for the model used in Theorem \ref{probabilistic version of main result}, as this model seems to be simpler to analyse.

Finally, we close the discussion of open problems with the following meta-problem. Theorems \ref{probabilistic version of main result} and \ref{deterministic version of main result} are quite similar in spirit in that the sequence $a$ is in both cases obtained by slightly perturbing an arithmetic progression by a suitably chosen $\{-1,0,1\}$-valued sequence $\varepsilon$. This prompts us to ask the following (somewhat vague and open-ended) question.

\begin{question}
\label{common generalisation of probabilistic and deterministic constructions}
Is there a common generalisation of Theorems \ref{probabilistic version of main result} and \ref{deterministic version of main result}? In particular, is there a notion of regularity for the sequence $\varepsilon$ that guarantees $a$ to have many consecutive sums and is satisfied by the instances in Theorems \ref{probabilistic version of main result} (with high probability) and \ref{deterministic version of main result}?
\end{question}
 
\noindent\textbf{Acknowledgements.} The author is grateful to Rudi Mrazović for his continued guidance and encouragement, as well as for useful discussions. He would also like to thank Thomas Bloom and Jakub Konieczny for helpful comments.

\end{document}